\documentclass[12pt]{article}
\usepackage{amsmath, amsthm, amssymb, epsfig, psfrag, amsfonts, latexsym}

\newtheorem{theorem}{Theorem}
\newtheorem{thm}[theorem]{Theorem}

\newtheorem{prop}[theorem]{Proposition}
\newtheorem{lemma}[theorem]{Lemma}

\newtheorem{Question}[theorem]{Question}

\theoremstyle{remark}

\theoremstyle{definition}

\def\cal{\mathcal}

\def\Bbb{\mathbb}

\def\D{\Bbb{D}}

\def\N{\Bbb{N}}

\def\Area{\hbox{\rm Area}}

\def\Diam{\hbox{\rm Diam}\,}

\def\FL{\hbox{\rm FL}}
\def\F+L{\hbox{$\textup{F}\!_+\textup{L}$}}

\def\ssm{\smallsetminus}

\def\onto{{\kern3pt\to\kern-8pt\to\kern3pt}}

\def\<{\langle}
\def\>{\rangle}
\def\|{{\ |\ }}

 \def\TT{\cal T}

\newcommand{\set}[1]{\left\{#1\right\}}

\renewcommand{\ni}{\noindent}
\renewcommand{\ss}{\smallskip}
\newcommand{\ms}{\medskip}

\def\*{^{\ast}}

\title{The absence of efficient  dual pairs \\ of spanning trees in planar graphs}

\author{T.R.Riley and W.P.Thurston
\thanks{The authors gratefully acknowledge support from NSF grants DMS--0540830 and DMS--0513436.} \\  
\small Mathematics Department, 310 Malott Hall, Cornell University, Ithaca, NY 14853-4201, USA  \\[-0.8ex]  \small \texttt{tim.riley@math.cornell.edu},  \texttt{wpt@math.cornell.edu}   
}

\date{November 2005; revised August 2006}

\begin{document}
\maketitle

\begin{abstract}
\ni A spanning tree $T$ in a finite planar connected graph $G$ determines a dual spanning tree $T\*$ in the dual graph $G\*$ such that $T$ and $T\*$ do not intersect.  We show that it is not always possible to find  $T$ in $G$ such that the diameters of $T$ and $T\*$ are both within a uniform multiplicative constant (independent of $G$) of the diameters of their ambient graphs.  

\ss \ni \small  2000 Mathematics Subject Classification: 05C10, 05C12, 20F06, 57M15 

\end{abstract}

\section{Introduction}\label{intro}

\begin{figure}[ht]
\psfrag{D}{$T^{\ast}$}
\psfrag{T}{$T$}
\psfrag{G}{$G$}
\centerline{\epsfig{file=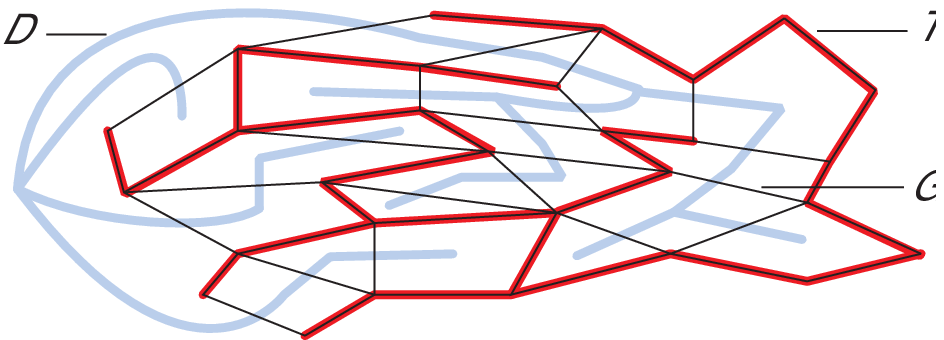}} 
\caption{Dual spanning trees.} \label{Dual pair} 
\end{figure}

Suppose $G$ is a finite connected undirected graph (or multigraph) embedded in the plane. Given a spanning tree $T$ in $G$, define $T\*$ to be the spanning tree in the dual graph $G\*$ whose edges are those dual to edges in $G \ssm T$.   Figure~\ref{Dual pair} gives an example.  

The length of a walk in a graph is the number of edges it contains and the distance between two vertices is the length of the shortest walk between them.  The diameter $\Diam G$ of a finite connected graph $G$ is the maximum distance between pairs of vertices of $G$.  

Motivated by issues arising in Geometric Group Theory concerning the geometry of van~Kampen diagrams, Gersten \& Riley  asked \cite{GR4}:
\begin{Question} \label{qn}
Does there exists $C>0$ such that if $G$ is a finite connected planar \textup{(}multi-\textup{)} graph then there is a maximal tree $T$ in $G$ with
\begin{eqnarray*}
\Diam T & \leq &  C \, \Diam G,  \text{ and} \\ 
\Diam T\* & \leq &  C \, \Diam G\*?
\end{eqnarray*}
\end{Question}
\ni They conjectured positive answers to a number of variants of this question with bounds imposed on the degrees of vertices in $G$ or $G\*$.  
We exhibit a family of graphs resolving these negatively.

\begin{thm} \label{bounds thm}
There are families $(G_n)_{n \in \mathbb{N}}$ of finite connected planar graphs such that all vertices in $G_n$ and $G\*_n$ have degree at most $6$, and there are constants $C_1,C_2 >0$ such that for all $n \in \N$ and all spanning trees $T$ in $G_n$,
\begin{eqnarray}
 \Diam G_n + \Diam G\*_n & \leq &  C_1 n, \, \text{ and} \label{diam bound} \\ 
 \Diam T + \Diam T\* & \geq &  C_2 n^2.  \label{lower bound}
 \end{eqnarray}
\end{thm}

Establishing (\ref{lower bound}) involves two key ideas.  The first is to regard $G_n$ as the 1-skeleton of a combinatorial 2-disc $\Delta_n$ and invoke a concept known as \emph{filling length}.  In the context of a simply connected metric space, Gromov~\cite{Gromov} defined the \emph{filling length} of a based loop $\gamma$ to be the infimal $L$ (assuming it exists) such that $\gamma$ can be contracted  through a family of based loops each of length at most $L$ to the constant loop (i.e.\ to the basepoint).    
We will use a combinatorial analogue of filling length from \cite{GR1} concerning \emph{shellings} of \emph{diagrams}. 

A \emph{diagram} $(\Delta,\star)$ is a finite planar contractible combinatorial 2-complex $\Delta$ equipped with a base vertex $\star$ on its boundary.  One can regard $\Delta$ as a finite planar multigraph $G$, the 1-skeleton of $\Delta$, with a 2-cell filling each face other than the \emph{outer} (i.e.\ unbounded) face.   Define the boundary walk of $\Delta$ based at $\star$ to be the  anti-clockwise  closed walk around the boundary of $\Delta$ that has origin $\star$ and follows the attaching map of the outer face.  The length of the boundary walk is the number of edges it contains (note that those in 1-dimensional portions of $\Delta$ are counted twice), or equivalently the degree of the vertex of $G\*$ dual to the outer face of $G$.

A \emph{shelling} of a diagram $\Delta = \Delta^0$ down to a vertex $\star$ on its boundary is a sequence $(\Delta^i)_{i=0}^m$ of diagrams in which $\Delta^m$ is the single vertex $\star$ and, for all $i$, we obtain $\Delta^{i+1}$ from $\Delta^i$ by one of the following two moves. 
\begin{itemize}
\setlength{\itemsep}{3pt} \setlength{\parsep}{3pt}
\item Remove a pendent edge and incident leaf $v \neq \star$. 
\item Remove an edge $e$ and the interior of a (closed) 2-cell $f$ where $e$ is in the boundaries of both $f$ and  $\Delta^i$.  
\end{itemize}
Each such move results in an elementary homotopy of the boundary walk: in the first case a backtracking pair of edges is removed, and in the second $e$ is replaced by the complementary portion of the walk around the boundary of $f$.  These moves ultimately achieve the contraction of the boundary walk of $\Delta$ down to the trivial walk at $\star$.   So we define the \emph{filling length} $\FL(\Delta, \star)$ of $(\Delta, \star)$ to be the minimal $L$ such that there is a shelling $(\Delta^i)_{i=0}^m$ of $\Delta$ in which for all $i$, the length of the boundary walk of $\Delta^i$ is at most $L$.  

Filling length will be useful to us because, given a diagram $(\Delta,\star)$ with $G$ the 1-skeleton of $\Delta$, the layout of a spanning tree $T$ in $G$ and the corresponding $T\*$ in $G\*$ can be made to dictate a shelling of $\Delta$ with filling length bounded above in terms of  $\Diam T + \Diam T\*$  (see Proposition~\ref{two trees}).  So a lower bound on the filling length of $(\Delta, \star)$ leads to a lower bound on $\Diam T + \Diam T\*$.

This brings us to the second key idea, which is to construct diagrams $(\Delta_n, \star)$ so as to contain a \emph{fattened tree} that forces the filling length of  $(\Delta_n, \star)$ to be suitably large.  In the context of Riemannian 2-discs this has been done by Frankel \& Katz in \cite{Frankel-Katz}, answering a question of Gromov; our $\Delta_n$ will essentially be combinatorial analogues of their metric discs.  To obtain $\Delta_n$ we first inductively define a family of trivalent trees $\TT_n$ by taking $\TT_0$ to be a lone edge, and $\TT_n$ to be three copies of $\TT_{n-1}$ with a leaf of each identified.  (We note that this does not determine $\TT_n$ uniquely.)  We then \emph{fatten} $\TT_n$ to a complex $A_n$ (see Figure~\ref{fat trees}) in which each of its edges becomes an $n \times n$ grid.  Finally, to obtain $\Delta_n$ we attach a combinatorial \emph{hyperbolic skirt} (a planar 2-complex $B_n$ that is topologically an annulus -- see Figure~\ref{hyp disc}) around the boundary of $A_n$ to reduce the diameter of its 1-skeleton to $\sim\! n$.  

Imagine inscribing $\TT_n$ in the plane, circling it with a loop, and then contracting that loop down to a point.  In the course of being contracted, the loop will intersect $\TT_n$.  In Lemma~\ref{many intersections} we show that however the loop contracts it must, at some time, meet at least $n+1$ distinct edges of $\TT_n$.  Envisage $A_n$ to be inscribed with a copy of $\TT_n$ as in Figure~\ref{fat trees}.  The lemma can be applied to the boundary walks  of the diagrams $\Delta^i_n$ of any shelling of $\Delta_n$ to learn that for some $i$ at least $n+1$ distinct edges of $\TT_n$ will be intersected; it then follows from the construction of $\Delta_n$ that at that time the length of the boundary walk is $\Omega(n^2)$. 

\ms
\ni \emph{Acknowledgement.}  Question~\ref{qn} was a topic of class discussion in a course taught by the second author at Cornell University in the Fall, 2005.  We are grateful to the members of the class, particularly John Hubbard and Greg Muller, for their contributions.  Additionally, we thank Andrew Casson, Genevieve Walsh and two anonymous referees for their comments on earlier versions of this article.

\section{Constructing the graphs $G_n$}

\begin{figure}[ht] 
\psfrag{p}{$p$}
\centerline{\epsfig{file=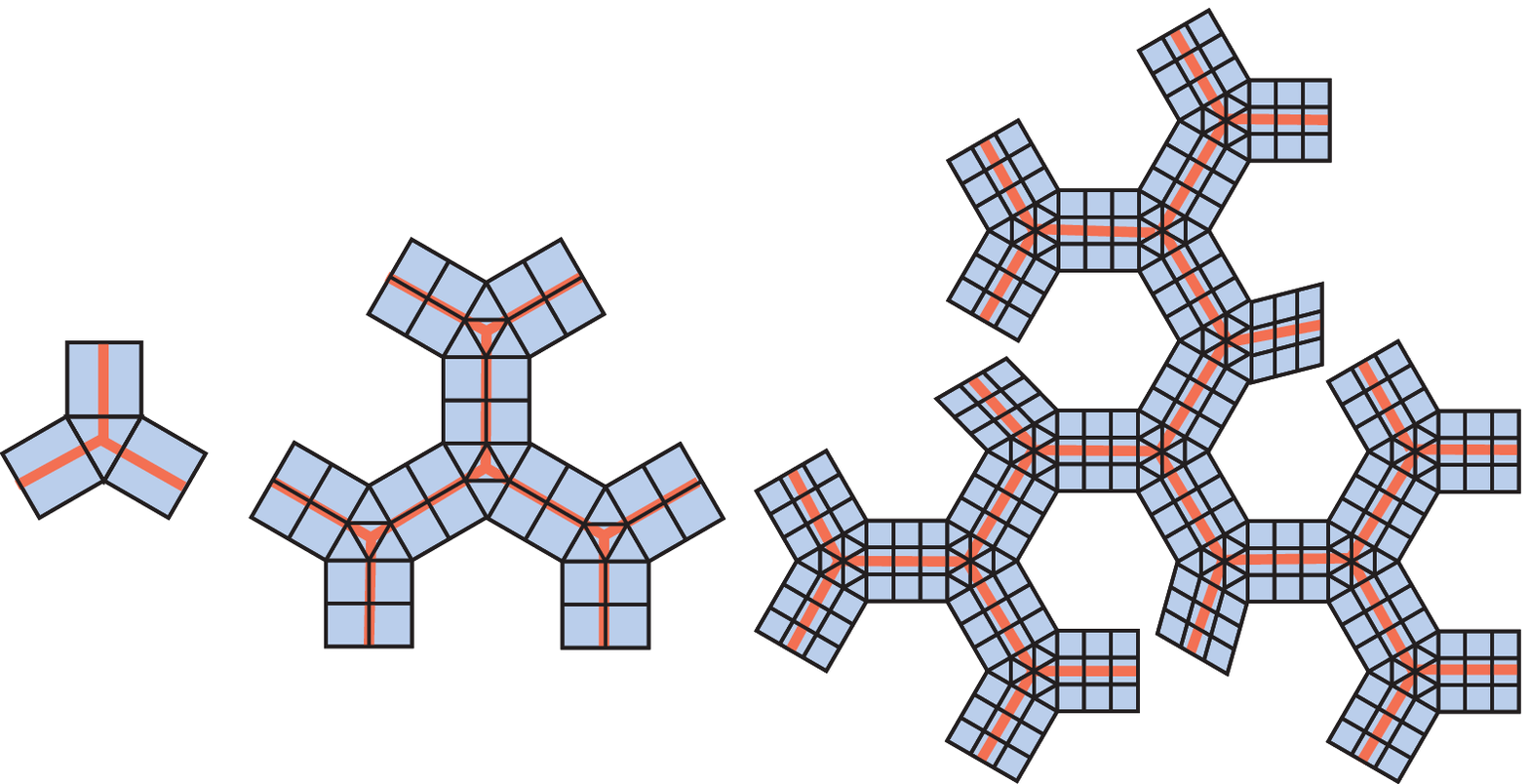}}
\caption{$A_1$, $A_2$ and $A_3$ inscribed with $\TT_1$, $\TT_2$ and $\TT_3$.} \label{fat trees} 
\end{figure}

Let $A_n$ be the family of diagrams (\emph{fattened trees}) obtained from $\TT_n$ (shown underlying) as illustrated in Figure~\ref{fat trees} by replacing edges by $n \times n$ grids and non-leaf vertices by tessellated triangles. 

\begin{figure}[ht] 
\psfrag{8}{$D_8$}
\psfrag{16}{$D_{16}$}
\psfrag{32}{$D_{32}$}
\psfrag{44}{$D_{44}$}
\centerline{\epsfig{file=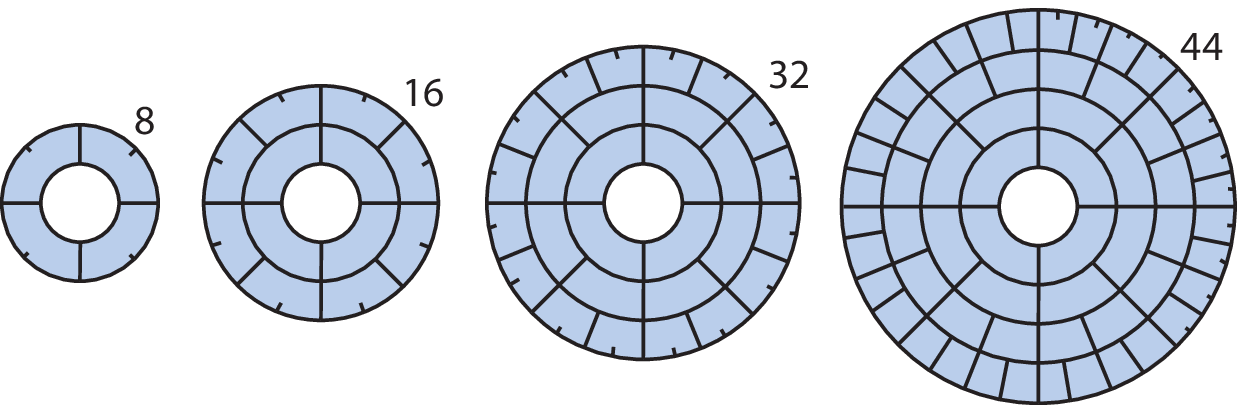}}
\caption{The annular 2-complexes $D_k$.} \label{hyp disc} 
\end{figure}

For $k=2^m$ with $m \geq 3$ define $D_k$ to be the planar combinatorial 2-complex that is topologically an annulus and is built out of $m-2$ concentric rings of pentagons as shown in Figure~\ref{hyp disc} for $m=3,4,5$.   For $m \geq 3$ and $2^{m-1}  < k \leq 2^m$, obtain $D_k$ from $D_{2^m}$ by inserting single edges in place of pairs of adjacent edges sharing a degree--two vertex until the total number of edges in the outer boundary cycle is reduced to $k$.  Figure~\ref{hyp disc} shows the example of $D_{44}$. 

The combinatorial length of the boundary circuit of $A_n$ is $p_n:=(5. 3^n + 3)n/2$.  For $n \geq 1$, define $B_n := D_{p_n}$, which plays the role of a \emph{hyperbolic skirt}: attach $A_n$ to $B_n$ by identifying the boundary of $A_n$ with the \emph{outer} boundary circuit of $B_n$ to give the planar combinatorial 2-disc $\Delta_n$.  Let $G_n$ be the 1-skeleton of $\Delta_n$.    

\section{Diameter estimates} \label{estimates}

We will now show that $(G_n)_{n \in \N}$ enjoys the properties listed in Theorem~\ref{bounds thm}.
By inspection, every vertex in $G_n$ and $G^{\ast}_n$ has degree at most $6$.  
Every vertex in $A_n$ is a distance at most $(n - 1)$ from the boundary, and one checks that the diameter of $B_n$ is at most a constant times $n$ since the number of concentric rings is $O(\log p_n)$.  Combined with similar considerations for the dual graphs this shows that there exists $C_1$ for which (\ref{diam bound}) holds.

For (\ref{lower bound}) we will use the following inequality  from \cite{GR3} on filling length.   (In fact, the definition of a shelling used in \cite{GR1, GR3} allows a third move, omitted from our the definition in Section~\ref{intro}, but that move is not needed here and plays no role in the proofs of the results cited in this article, namely Propositions~\ref{two trees} and \ref{log shelling prop}.)

\begin{prop}[Proposition~3.4, \cite{GR3}] \label{two trees}
Suppose $(\Delta,\star)$ is a diagram in which the degree of each 2-cell is at most $\lambda$.  If $T$ is a spanning tree in the 1-skeleton of $\Delta$ then 
\begin{eqnarray} \label{corrected bound}
\FL(\Delta,\star) \  \leq \ \Diam T  \, +  \, 2\, \lambda \, \Diam T\*  \, + \,  \ell(\partial\Delta), 
\end{eqnarray}
where $\ell(\partial\Delta)$ denotes the length of the boundary walk of $\Delta$. 
\end{prop}

We refer the reader to \cite{GR3} for a detailed proof, but will sketch the idea here. 
Regard the vertex of $T\*$ outside $\Delta$ as the root $r$ of $T\*$.  The embedding of $T\*$ in the plane defines a cyclic ordering on its leaves.  Define a $T\*$-gallery of $\Delta$ to be a subcomplex that is the union of the closed 2-cells of $\Delta$ that are dual to the vertices lying on a path in $T\*$ from $r$ to a leaf. 
The idea is that \emph{tunnelling} along paths of $T\*$ from $r$ to successive leaves, following their cyclic ordering, dictates a shelling $(\Delta^i)$ of $\Delta$ that establishes (\ref{corrected bound}):  when traversing an edge $e\*$ in such a path shell the edge $e$ dual to $e\*$ and the face dual to the terminal vertex of $e\*$;  en route, remove all pendant edges (with leaf vertices $\neq \star$) immediately they become available.    The boundary walks of the diagrams $\Delta^i$ are then each comprised of a path in $T$, trails in the 1-skeleta of two $T\*$-galleries of $\Delta^i$, and a portion of the boundary walk of $\Delta$.  Thus we get (\ref{corrected bound}).
\ms

For the following lemma and subsequent discussion it is convenient to regard $\TT_n$ as a \emph{disjoint} union of its edges; accordingly choose one edge in $\TT_n$ to include both of its end-vertices and all others to include exactly one end-vertex. 

\begin{lemma} \label{many intersections}
Suppose $\TT_n$ is embedded in a disc, which for convenience we take to be the unit disc in the complex plane.  Suppose $H : [0,1]^2 \to \D^2$ is  a continuous map \textup{(}a homotopy\textup{)} satisfying $H(0,t) = H(1,t) = 1$ for all $t$, and $H_0(s) = e^{2 \pi i s}$ and $H_1(s) = 1$ for all $s$,  where $H_t$ denotes the restriction of $H$ to $[0,1] \times \set{t}$.   Further, assume $H([0,1] \times [0,t]) \cap  H([0,1] \times [t,1]) = H([0,1] \times \set{t})$ for all $t$. Then $H_t$ meets at least $n+1$ edges in $\TT_n$ for some $t \in [0,1]$.
\end{lemma}

\begin{proof} 
The case $n=0$ is immediate.  For the induction step, express $\TT_n$ as the wedge $\bigwedge_{i=1}^3 \TT^i_{n-1}$ of three copies of $\TT_{n-1}$ at a vertex $v$.  Obtain $\hat{\TT}^i_{n-1}$ from $\TT^i_n$ by removing a small open neighbourhood of $v$.   
Let $t_i$ be such that $H_{t_i}$ meets at least $n$ edges of $\hat{\TT}^i_{n-1}$.  Renumbering if necessary,  we may assume $t_1 \leq t_2 \leq t_3$.  The condition that $H([0,1] \times [0,t]) \cap  H([0,1] \times [t,1]) = H([0,1] \times \set{t})$ for all $t$, ensures 
that if $t_1 \leq t \leq t_3$ and $H_t([0,1]) \cap (\TT^1_{n-1} \cup \TT^3_{n-1}) = \emptyset$ then points of $H_{t_1}([0,1]) \cap (\TT^1_{n-1} \cup \TT^3_{n-1})$ are in different path components than points $H_{t_3}([0,1]) \cap (\TT^1_{n-1} \cup \TT^3_{n-1})$ in $\D^2 \ssm H_{t_3}([0,1])$, but that is impossible as $\TT^1_{n-1} \cup \TT^3_{n-1}$ is path connected.  We deduce, in particular, that $H_{t_2}$ intersects $\TT^1_{n-1} \cup \TT^3_{n-1}$, and so meets at least $n+1$ edges of $\TT_n$.      
\end{proof}

We can now establish (\ref{lower bound}).  Choose any vertex on the boundary of $\Delta_n$ to serve as the base vertex  $\star$.  Envision the subdiagram $A_n$ of $\Delta_n$ to be inscribed with $\TT_n$ as  in Figure~\ref{fat trees}.  The diagrams $\Delta^i_n$ of a shelling of $(\Delta_n, \star)$ are subcomplexes whose boundary walks define concentric loops ultimately contracting to $\star$.  Interpolating suitably between these loops produces a homotopy in which the boundary walk of $\Delta_n$ is contracted to the constant loop at $\star$ through a family of loops $H_t$.  So by Lemma~\ref{many intersections} there exists $t$ such that $H_t$  meets $n+1$ edges of $\TT_n$ and it follows that there exists $i$ such that the boundary walk of $\Delta_n^i$ meets $n+1$ edges of $\TT_n$.  But any path in the 1-skeleton of $\Delta_n$ meeting four distinct edges of $\TT_n$ has combinatorial length at least $n$.  So the length of the boundary walk of $\Delta_n^i$ is at least $n \lfloor n / 3 \rfloor$.  
Deduce that $\FL(\Delta_n, \star) \geq n \lfloor n / 3 \rfloor$ and therefore, by Proposition~\ref{two trees}, there exists $C_2>0$ such that (\ref{lower bound}) holds.

\section{Two concluding remarks}

We note that Proposition~3.3 in \cite{GR3} exhibits another family of diagrams in which filling length outgrows 1-skeleton diameter. However, filling length does not outgrow the diameter of the dual in these examples.    

Finally, we mention that our family of diagrams $\Delta_n$ exhibits the most radical divergence possible between filling length, diameter and dual diameter in the sense of the following result.

\begin{prop} \label{log shelling prop}
Given $\lambda>0$, there exists $C = C(\lambda)$ such that if $(\Delta,\star)$ is a diagram in which the degree of each 2-cell is at most $\lambda$ then
$$\FL(\Delta, \star),  \  \leq  \ C (\Diam G)(\Diam G\*),$$ 
where $G$ is the 1-skeleton of $\Delta$.  
\end{prop}

This follows from an argument of \cite{GR1} which we will only briefly outline here.  Take a \emph{geodesic} spanning tree $T$ in $G$ based at $\star$ -- that is, a spanning tree such that for all vertices $v$ in $G$, the distance from $v$ to $\star$ in $T$ is the same as in $G$.  Note that $\Diam T \leq 2 \Diam G$.  Let the vertex $r$ of $G\*$ that is outside $\Delta$ be the root of $T\*$.  By \emph{subtrees suspended from a vertex $v$} in $T\*$ we mean the closures of the connected components of $T\* \ssm \set{v}$ that do not contain $r$.  Describe a vertex as \emph{branching} when there is more than one  subtree suspended from it.   A vertex \emph{below} $v$ is any vertex of any subtree suspended from $v$.  Define the \emph{weight} of a tree to be the number of vertices it contains that have degree at least three.  

Consider \emph{tunnelling through} $\Delta$ along the walk in $T\*$ that starts at $r$, first proceeds to the nearest leaf or branching vertex (possibly $r$ itself) and then continues  according to the following rules from its current vertex $v$.    
\begin{itemize}
\item If $v$ is a branching vertex then of the as-yet-unentered subtrees suspended from $v$, choose one of least weight and proceed to the nearest leaf or branching vertex ($\neq v$) therein. 
\item If $v$ is a leaf return to the most recently visited branching vertex attached to which there remain as-yet-unentered suspended subtrees of $T\*$.  
\end{itemize}
The walk is complete when every edge in $T\*$ has been traversed.
This walk dictates the following shelling of $\Delta$ (termed \emph{logarithmic shelling} in \cite{GR1}): when traversing an edge $e\*$ for the first time, remove the dual edge $e$ and the face dual to the terminal vertex of $e\*$, and immediately any pendant edge (with leaf vertex not $\star$)  appears, remove it.   

The lengths of the boundary walks of the diagrams encountered in this shelling are at most a constant (depending on $\lambda$) times $(\Diam T) \log (1+\Area \, \Delta)$, where $\Area \, \Delta$ denotes the number of 2-cells in $\Delta$.  As $\Area \, \Delta  \leq \lambda^{\mbox{$\Diam G\*$}}$ and $\Diam T \leq 2 \Diam G$, the result follows.  

\bibliographystyle{plain}
\bibliography{bibli}

\def\cprime{$'$}
\begin{thebibliography}{1}

\bibitem{Frankel-Katz}
S.~Frankel and M.~Katz.
\newblock The {M}orse landscape of a {R}iemannian disc.
\newblock {\em Ann. Inst. Fourier, Grenoble}, 43(2):503--507, 1993.

\bibitem{GR1}
S.~M. Gersten and T.~R. Riley.
\newblock Filling length in finitely presentable groups.
\newblock {\em Geom. Dedicata}, 92:41--58, 2002.

\bibitem{GR4}
S.~M. Gersten and T.~R. Riley.
\newblock Some duality conjectures for finite graphs and their group theoretic
  consequences.
\newblock {\em Proc.\ Edin.\ Math.\ Soc.}, 48(2):389--421, 2005.

\bibitem{GR3}
S.~M. Gersten and T.~R. Riley.
\newblock The gallery length filling function and a geometric inequality for
  filling length.
\newblock {\em Proc. London Math. Soc.}, 92(3):601--623, 2006.

\bibitem{Gromov}
M.~Gromov.
\newblock Asymptotic invariants of infinite groups.
\newblock In G.~Niblo and M.~Roller, editors, {\em Geometric group theory II},
  number 182 in LMS lecture notes. Camb. Univ. Press, 1993.

\end{thebibliography}

\end{document}